\documentclass[a4paper,11pt,reqno]{article}
\usepackage{amsmath, amsfonts, amssymb, amsthm}
\usepackage{mathrsfs}
\usepackage{amsmath}
\usepackage{amssymb,color}
\usepackage{graphicx}
\usepackage{amssymb}
\usepackage{authblk}
\usepackage[numbers,sort&compress]{natbib}

\usepackage[colorlinks=true]{hyperref}
\hypersetup{urlcolor=blue, citecolor=blue}

\par

\def \R {{\mathbb{R}}}
\numberwithin{equation}{section}

\begin{document}

\title{On the rotation-two-component $b$-family shallow-water system: derivation and wave-breaking analysis}

\author[1,2]{Xingxing Liu\thanks{E-mail: liuxxmaths@cumt.edu.cn}}

\affil[1]{School of Mathematics, China University of Mining and Technology,\newline
 Xuzhou, Jiangsu 221116, China}
\affil[2]{Jiangsu Center for Applied Mathematics (CUMT),\newline
 Xuzhou, Jiangsu 221116, China}

\date{}
\maketitle

\begin{abstract}
In this paper, we derive the rotation-two-component $b$-family shallow-water wave system from the Green--Naghdi type equations incorporating the effects of the Coriolis force and constant vorticity.
The system admits, among its special cases, the two-component Degasperis--Procesi (2DP) system as a reduction, thereby confirming its genuine role as a rotational extension of the
$b$-family.
We then focus on a specific reduction, the rotation-two-component Degasperis-Procesi system, whose particular interest lies in the fact that it combines two independent sources of difficulty:
unlike the two-component Camassa-Holm type system, it possesses no conserved $H^1\times L^2$-type quantity, and in contrast to the standard 2DP system, it contains a higher-order nonlinear coupling term induced by the Coriolis force.
To overcome this  double difficulty, we introduce an auxiliary variable and establish a local $L^\infty$-estimate for $u$ without using any conservation law, based on which we derive two distinct blow-up criteria. When the Coriolis force vanishes, our results recover known criteria for the 2DP system and also provide new criteria for both the 2DP system and Degasperis-Procesi equation. The results extend naturally to the periodic setting and to relevant special cases of the rotational $b$-family system.

\noindent Mathematics Subject Classification: 35Q35; 35B44; 35G25; 76B15
\smallskip\par
\noindent \textit{Keywords}: Rotation-two-component $b$-family system; Rotation-two-component Degasperis--Procesi system; Coriolis force; Wave-breaking

\end{abstract}

\section{Introduction}
\newtheorem{theorem1}{Theorem}[section]
\newtheorem{lemma1}{Lemma}[section]
\newtheorem {remark1}{Remark}[section]
\newtheorem {definition1}{Definition}[section]
\newtheorem{corollary1}{Corollary}[section]
\newtheorem{proposition1}{Proposition}[section]
\par

The Camassa-Holm (CH) equation is a nonlinear dispersive wave equation that arises in shallow-water theory. It is derived by approximating directly the Hamiltonian for Euler equations in the shallow water regime \cite{CH}.
 Subsequently, it was established as a valid approximation to the Euler equations for modeling the unidirectional propagation of waves at the free surface of shallow water under the influence of gravity \cite{CL,Johnson}.
The CH equation is completely integrable, admits peaked solitary wave solutions (peakons), and exhibits wave-breaking \cite{CE,CL1,CS}. These remarkable properties have made it an important model in shallow-water wave theory and have motivated a number of generalisations.

On the one hand, the single component $b$-family equation \cite{HS,HS1}
\begin{equation}\label{bfamily}
m_t + u m_x + b u_x m = 0,\qquad m = u - u_{xx},
\end{equation}
extends the CH equation by introducing an extra parameter $b$, which regulates the balance between convection $u m_x$ and stretching $b u_x m$.  For $b=2$ it reduces to the CH equation, while $b=3$ recovers the Degasperis-Procesi (DP) equation \cite{DP}.

On the other hand, the CH equation has been extended to two-component systems by coupling an additional variable $\rho$ that is related to the free surface deviation.
Two main derivation routes have emerged: one starts from the Green-Naghdi (GN) equations  \cite{ASL,GN} and performs an asymptotic expansion \cite{CI}; the other begins directly from the Euler equations with constant vorticity \cite{Ivanov}.

It is known that
geophysical water waves in equatorial regions are strongly influenced by the Earth's rotation. The Coriolis force plays a key role in large-scale ocean dynamics, such as the generation of azimuthal waves along the equator \cite{CJ1,CJ2}. At the same time, the equatorial undercurrent (EUC) is a prominent feature of the Pacific Ocean \cite{CI1}. Interestingly, the effect of the Coriolis force on the propagation of waves across the equatorial Pacific is relatively weak because the latitudinal variation of the EUC is small. Nevertheless, the rich physical phenomena observed over a few degrees of latitude near the equator motivate the development of shallow-water models that can capture both the Coriolis effect and the presence of the EUC \cite{HM,HYZ}.

To describe such flows, the GN-type equations have been progressively generalised: first a pure-Coriolis GN equations \cite{GLL}, then a pure-shear GN equations \cite{ZL}, and most recently a full GN-type equations \cite{HYZ} that incorporates both the Coriolis parameter $\Omega$ and the constant vorticity $A$. In the present paper we follow the GN route, but we take this full GN-type equations as our starting point.

Two small dimensionless parameters govern the shallow-water approximation: the amplitude parameter $\varepsilon = a/h_0$ and the shallowness parameter $\mu = h_0^2/\lambda^2$, where $h_0$ is the mean water depth, $a$ the typical wave amplitude and $\lambda$ the typical wavelength. We work under the Boussinesq scaling $\mu\ll1$, $\varepsilon = O(\mu)$, which corresponds to weakly nonlinear long waves. The alternative CH scaling $\varepsilon = O(\sqrt{\mu})$ \cite{CL} would lead to stronger nonlinearity and wave-breaking, but lies beyond the scope of the present work. Expanding the full GN-type equations to order $O(\mu)$ yields a simplified system in which the equation for $\eta$ lacks the dispersive term $\frac{\mu}{6}u_{xxx}$ that appears in the asymptotic systems directly derived from the Euler equations under the same scaling (e.g. \cite{EHKL,FGL}). This absence makes the introduction of an auxiliary variable considerably cleaner.

By performing a systematic asymptotic expansion up to $O(\varepsilon^2,\varepsilon\mu,\mu^2)$ and a Galilean transformation $x\to x-At$, $t\to t$, we obtain, through a series of algebraic manipulations and cancellations of higher-order terms, the following rotation-two-component $b$-family (R2$b$-family) shallow-water wave system
\begin{equation}\label{systemsus}
\left\{
 \begin{aligned}
&m_t-Au_x+\sigma(bu_xm+um_x)+(b+1)(1-\sigma)uu_x\\
&\quad\quad\quad\quad\quad\quad\quad\quad\quad\quad+(1-2\Omega A)\rho\rho_x-2\Omega\rho\big((b-1)\rho u_x+u\rho_x\big)=0,\\
&\rho_t+u\rho_x+(b-1)\rho u_x=0,
\end{aligned}
\right.
\end{equation}
where $m=u-u_{xx}$. The variable $u(t,x)$ represents the horizontal velocity of the fluid, and $\rho(t,x)$
is related to the free surface deviation from equilibrium, with the boundary conditions $u\to0$, $\rho\to1$ as $|x|\to\infty$.
Here $b\neq1$ is a real parameter, the constant $A$ denotes the constant vorticity, characterising a linear underlying shear flow, and the parameter $\Omega$ is the constant rotational frequency due to the Coriolis effect.
The parameter $\sigma$, inherited from hyperelastic rod models \cite{D,DH} and the generalized two-component Camassa-Holm (2CH) system \cite{CLiu}, controls the balance between nonlinear steepening and amplification in fluid convection due to stretching.

The R2$b$-family system (\ref{systemsus}) unifies several known models, as detailed below.
\begin{itemize}
\item The rotation-two-component CH system (R2CH) \cite{FGL}.  Setting $b=2$ in (\ref{systemsus}) gives
\begin{equation}\label{systemsb2}
\left\{
 \begin{aligned}
&m_t-Au_x+\sigma(2u_xm+um_x)+3(1-\sigma)uu_x+(1-2\Omega A)\rho\rho_x-2\Omega \rho(\rho u)_x=0,\\
&\rho_t+(\rho u)_x=0.
\end{aligned}
\right.
\end{equation}
Note that this system (\ref{systemsb2}) was derived from the Euler equations, whereas our derivation starts from the GN-type equations. The two approaches lead to the same model.

\item The constant vorticity two-component $b$-family system \cite{EHKL}.  Taking $\Omega=0$ and $\sigma=1$ in (\ref{systemsus}) yields
\begin{equation}\label{systemsconstant2}
\left\{
 \begin{aligned}
&m_t - A u_x + (b u_x m + u m_x) + \rho\rho_x = 0,\\
 &\rho_t + u\rho_x + (b-1)\rho u_x = 0.
\end{aligned}
\right.
\end{equation}
Like the R2CH system (\ref{systemsb2}), this system (\ref{systemsconstant2}) was also derived from the Euler equations and, in contrast to our system (\ref{systemsus}), does not feature the Coriolis force.
Obviously, it can reduce to 2CH with $b=2$, and two-component DP (2DP) with $b=3$ \cite{P}, respectively.
\item The generalized 2CH system \cite{CLiu}. When $b=2$ and $\Omega=0$, the system (\ref{systemsus}) reduces to
\begin{equation*}\label{systemsgeneralized2}
\left\{
 \begin{aligned}
&m_t - A u_x + \sigma(2u_x m + u m_x) + 3(1-\sigma)u u_x + \rho\rho_x = 0,\\
&\rho_t +(\rho u)_x = 0.
\end{aligned}
\right.
\end{equation*}
When $\sigma=1$, this system further reduces to the constant vorticity 2CH system studied by Ivanov \cite{Ivanov}; when additionally $A=0$, we recover the classical 2CH system of Constantin and Ivanov \cite{CI}.

\item The single component $b$-family equation \cite{HS,HS1}. In the limit $\rho\equiv0$, the first equation of (\ref{systemsus}) becomes
\[
m_t - A u_x + \sigma(b u_x m + u m_x) + (b+1)(1-\sigma)u u_x = 0,
\]
which, taking $\sigma=1$ and $A=0$, recovers the standard $b$-family equation (\ref{bfamily}).
\end{itemize}

It is also instructive to compare our R2$b$-family system (\ref{systemsus}) with the following rotational $b$-family (R-$b$-family) system proposed by Zhu and Wang \cite{ZW}. Their system reads
\begin{equation}\label{systemsoriginal}
\left\{
 \begin{aligned}
&m_t-Au_x+\sigma(bu_xm+um_x)+(b+1)(1-\sigma)uu_x
+(1-2\Omega A)\rho\rho_x-2\Omega \rho(\rho u)_x=0,\\
&\rho_t+(\rho u)_x=0.
\end{aligned}
\right.
\end{equation}
In the R-$b$-family system (\ref{systemsoriginal}), the coupling term $\rho(\rho u)_x$ in the first equation does not involve the parameter $b$
and the second equation lacks the stretching term $(b-1)\rho u_x$. Consequently, when $b=3, \Omega=0$ their system does not reduce to the 2DP system. In contrast, our system (\ref{systemsus}) has fully $b$-dependent couplings and therefore represents the genuine rotational extension of the $b$-family equation. Moreover, while the system (\ref{systemsoriginal}) was derived from the Euler equations, our derivation starts from the GN-type equations, providing a different asymptotic framework.

Therefore, the R2$b$-family system (\ref{systemsus}) unifies several known models within a single framework.
Upon examining its various parametric special cases,
we find that when $b=3,\; A=0,\; \sigma=1$, the system (\ref{systemsus}) reduces to the following rotation-two-component DP (R2DP) system, which takes the form
\begin{equation}\label{systems2}
\left\{
 \begin{aligned}
&m_t+um_x+3u_xm+\rho\rho_x-2\Omega \rho(2\rho u_x+u \rho_x)=0,\\
&\rho_t+2\rho u_x+u \rho_x=0.
\end{aligned}
\right.
\end{equation}
To the best of our knowledge, the R2DP system (\ref{systems2}) has not been studied in the literature, in contrast to the R2CH system, which has been extensively investigated.

This new system (\ref{systems2}) distinguishes itself from both the 2DP system studied in \cite{YY} and the R2CH system studied in \cite{FGL} by
combining two independent sources of mathematical difficulty. For the R2CH system, the presence of a conserved
$H^1\times L^2$-type conserved quantity provides crucial control over the nonlocal terms arising from the Coriolis force, making the auxiliary variable method in \cite{FGL} feasible.
For the 2DP system without rotation ($i.e.$, $\Omega=0$), although no such conservation law exists, the absence of the Coriolis-induced coupling term $\rho(2\rho u_x+u\rho_x)$
allows for certain structural identities that lead to an $L^\infty$-estimate for $u$ \cite{YY}.
In contrast, the R2DP system (\ref{systems2}) suffers from both difficulties simultaneously: the conservation law is absent, and the higher-order coupling term destroys the structural identities used in the 2DP case.
Consequently, our analysis must overcome a genuinely harder problem, and the local $L^\infty$-estimate for $u$
established in Lemma \ref{Lem3.3} is obtained without relying on either conservation laws or the particular structure of the 2DP system.

To overcome this double difficulty, we follow the idea of introducing an auxiliary variable as in \cite{FGL} and define
\[
\omega = u +\Omega\,p\ast \rho^2,\qquad p(x)=\frac12 e^{-|x|},
\]
and then derive the evolution equation for $\omega, \omega_x$ using system (\ref{systems2}). However, the situation here is very different from that in \cite{FGL}: due to the absence of a conservation law, we cannot directly control the boundedness of the terms in the $\omega$-equation as in the R2CH case.
To address this, we establish a local a priori estimate for the $L^\infty$-norm of $u$ for the R2DP system without using any conservation law (see Lemma \ref{Lem3.3}), thereby gaining the necessary control over the solution on a finite time interval.

Based on this a priori estimate, we derive two blow-up criteria for the R2DP system (see Theorems \ref{Th3.4} and \ref{Th3.5}). These two criteria differ in their mathematical nature: Theorem \ref{Th3.4} tracks the evolution of the $\omega_x$, which can be viewed as a generalization of the classical Riccati differential inequality approach; Theorem \ref{Th3.5} focuses on the dynamics of $\omega\pm\omega_x$ along characteristics, reflecting the use of the "local structure" of solutions. This latter idea originated from the work of Brandolese \cite{B}
on CH equation. With the $L^\infty$ a priori estimate for $u$ at hand, we are able to handle the extra nonlinearity induced by the Coriolis force and establish the corresponding blow-up criteria.

It is worth pointing out that
when the Coriolis force vanishes ($\Omega=0$), Theorem \ref{Th3.4} reduces to the classical blow-up criterion for the 2DP system given in \cite{YY}, while Theorem \ref{Th3.5} provides a new blow-up criterion for the 2DP system.
In the single-component limit $\rho\equiv 0$, both theorems yield corresponding blow-up criteria for the DP equation; in particular, Theorem \ref{Th3.5} gives a new blow-up result for the DP equation.
Moreover, as noted in Remarks \ref{Re3.2}-\ref{Re3.3}, the blow-up criteria established in this paper remain valid in the periodic setting and can also be extended to the relevant special cases of the R-$b$-family system.

The rest of this paper is organized as follows. In Section \ref{sec2}, we derive in detail the R2$b$-family system (\ref{systemsus}) and discuss its reductions when the parameters take special values. In Section \ref{sec3}, we first give the local well-posedness result, then present the precise blow-up scenarios; based on these, we introduce the auxiliary variable, establish the $L^\infty$ a priori estimate for $u$, and finally prove the two blow-up theorems together with their relevant remarks.

\section{Derivation of the model system}\label{sec2}
\newtheorem{theorem2}{Theorem}[section]
\newtheorem{lemma2}{Lemma}[section]
\newtheorem {remark2}{Remark}[section]
\newtheorem {definition2}{Definition}[section]
\newtheorem{corollary2}{Corollary}[section]
\par
In this section, we present a derivation of the R2$b$-family system (\ref{systemsus}) from the following GN-type equations with effects of the Coriolis force and equatorial undercurrent
\begin {equation}
\left\{
\begin{aligned}\label{rGNwithshear}
&{{\eta _t} + A(1 + \varepsilon \eta ){\eta _x} + {{\big((1 + \varepsilon \eta )u\big)}_x} = 0}, \\
&{u_t} + \varepsilon u{u_x} + {\eta _x}+2\Omega \eta _t
= \frac{\mu }{{3(1 + \varepsilon \eta )}}{{\big({{(1 + \varepsilon \eta )}^3}({u_{xt}} + \varepsilon u{u_{xx}} - \varepsilon u_x^2)\big)}_x}\\
&\quad \quad \quad \quad \quad \quad \quad \quad\quad \quad \quad  +\frac{A\mu}{9}\big(4((1 + \varepsilon \eta )^3u_{xx})_x-(1 + \varepsilon \eta )^3u_{xxx}\big)
+O({\mu ^2}).
\end{aligned}
\right.
\end{equation}
This GN-type equations (\ref{rGNwithshear}), recently derived in \cite{HYZ} from the Euler equations under the $f$-plane approximation,
incorporates both the Coriolis parameter $\Omega$ and the constant vorticity $A$ (generalizing earlier models with only $\Omega$ \cite{GLL} or only $A$ \cite{ZL}). It serves as the starting point for our asymptotic derivation.

We now perform the asymptotic expansion under the Boussinesq scaling
\begin{equation}\label{scaling}
\mu\ll1,\quad\quad \varepsilon=O(\mu),
\end{equation}
which corresponds to the weakly nonlinear long-wave regime.
Although $\varepsilon$ and $\mu$ are independent parameters, imposing $\varepsilon=O(\mu)$ balances nonlinearity and dispersion at leading order, leading to classical equations such as KdV or BBM \cite{CJ}. (For completeness, we note that the alternative CH scaling $\varepsilon = O(\sqrt{\mu})$ \cite{CL} would lead to stronger nonlinearity and wave-breaking phenomena, but it is not needed here because our derivation already captures the essential physics within the Boussinesq regime.)

Under the Boussinesq scaling (\ref{scaling}), we expand the second equation of (\ref{rGNwithshear}) in powers of $\mu$ and retain terms up to
 $O(\mu)$ (while neglecting higher-order terms of order $O(\varepsilon\mu,\mu^2)$). The equations (\ref{rGNwithshear}) then reduce to the following simplified system
\begin {equation}
\left\{
\begin{aligned}\label{jianhuarGNwithshear}
&{{\eta _t} +A\eta _x+ {{\big((1 + \varepsilon \eta )u+\varepsilon\frac{A}{2}\eta^2\big)}_x} = 0}, \\
&(u-\frac{\mu }{3}u_{xx})_t+ \varepsilon u{u_x}-\frac{A\mu}{3}u_{xxx}+ {\eta _x}+2\Omega \eta _t
=O(\varepsilon\mu,\mu^2).
\end{aligned}
\right.
\end{equation}
Notably, compared with the corresponding asymptotic system obtained in \cite{EHKL,FGL} under the same Boussinesq scaling, the dispersive term $\frac{\mu}{6}u_{xxx}$ is absent from the first equation of (\ref{jianhuarGNwithshear}). As will be seen, this simplifies the subsequent derivation.

We now linearise (\ref{jianhuarGNwithshear}) by neglecting all $\varepsilon$ and $\mu$ terms. This yields
\begin{equation*}
\left\{
\begin{aligned}
&{\eta _t}+A\eta_x+u_{x}=0,\\
&{u_{t}}+{\eta _x}+2\Omega\eta_t=0,
\end{aligned} \right.
\end{equation*}
giving us
\begin{equation}\label{linear}
\left\{
\begin{aligned}
&{\eta _{tt}}  - {\eta _{xx}}+(A-2\Omega)\eta_{tx} = 0,\\
&{u_{tt}} - {u_{xx}}+(A-2\Omega)u_{tx}  = 0.
\end{aligned}
\right.
\end{equation}
It is known that the equations (\ref{linear}) admit travelling wave solutions of the form $\eta = \eta(x-ct)$ and $u(x-ct)$.
Substituting into (\ref{linear}) gives the dispersion relation $c^2-(A-2\Omega)c-1=0$, and we select the right-running wave by taking the positive root
$c=\frac{A-2\Omega+\sqrt{(A-2\Omega)^{2}+4}}{2}$. Consequently, we obtain the following useful linear relations
\begin{equation}\label{relationueta}
\eta-\frac{1}{c-A}u=O(\varepsilon ,\mu ),\quad\eta_t +\frac{c}{c-A}u_x= O(\varepsilon ,\mu ),\quad
{u_{t}} +(1-2\Omega c){\eta _x} = O(\varepsilon ,\mu ),
\end{equation}
where we require $c\neq A$. Indeed, the case $c=A$ leads to $1-2\Omega A=0$ and forces $u$ to be constant, which we exclude.

To derive our R2$b$-family system (\ref{systemsus}), we introduce an auxiliary variable of the form
\begin{equation}\label{varrho}
\varrho=1+k_1\varepsilon\eta+k_2\varepsilon^2 \eta^2,
\end{equation}
where $k_1$ and $k_2$ are constants to be determined later. Solving (\ref{varrho}) for $\eta$ and its derivatives, we obtain
\begin{equation*}\label{etatx}
\left\{
\begin{aligned}
&\eta_t=\frac{1}{\varepsilon k_1}\varrho_t-\varepsilon\frac{2k_2}{k_1}\eta\eta_t,\\
&\eta_x=\frac{1}{\varepsilon k_1}\varrho_x-\varepsilon\frac{2k_2}{k_1}\eta\eta_x.
\end{aligned}
\right.
\end{equation*}
Using the linear relations (\ref{relationueta}) and keeping terms of order $O(\varepsilon,\mu)$, we find
\begin{equation}\label{etatx}
\left\{
\begin{aligned}
&\eta_t=\frac{1}{\varepsilon k_1}\varrho_t+\varepsilon\frac{2k_2}{k_1}\frac{c}{c-A}\eta u_x+O(\varepsilon^2,\varepsilon\mu),\\
&\eta_x=\frac{1}{\varepsilon k_1}\varrho_x-\varepsilon\frac{2k_2}{k_1}\frac{1}{c-A}\eta u_x+O(\varepsilon^2,\varepsilon\mu).
\end{aligned}
\right.
\end{equation}
We now substitute (\ref{etatx}) into the first equation of (\ref{jianhuarGNwithshear}). A straightforward calculation gives
\begin{eqnarray}\label{varrho1}
\frac{1}{\varepsilon k_1}(\varrho_t+A\varrho_x)+\varepsilon\frac{2k_2}{k_1}\eta u_x+{{\big((1 + \varepsilon \eta )u+\varepsilon\frac{A}{2}\eta^2\big)}_x}=O(\varepsilon^2,\varepsilon\mu).
\end{eqnarray}
Next, we eliminate $\eta$ in favour of $u$ using the relations (\ref{relationueta}). After neglecting terms of order
$O(\varepsilon^2,\varepsilon\mu)$, the equation (\ref{varrho1}) becomes
\begin{eqnarray}\label{varrho2}
\frac{1}{\varepsilon k_1}(\varrho_t+A\varrho_x)+\Big(\big(1+\varepsilon(1+\frac{1}{2}\frac{A}{c-A}+\frac{k_2}{k_1})\eta\big)u\Big)_x=0.
\end{eqnarray}
At this stage we are free to choose the constants $k_1, k_2$ so that the expression inside the parentheses simplifies.
We impose the constraint
\begin{eqnarray}\label{relation1}
1+\frac{1}{2}\frac{A}{c-A}+\frac{k_2}{k_1}=k_1.
\end{eqnarray}
With this choice, and recalling the definition (\ref{varrho}) of $\varrho$, the equation (\ref{varrho2}) reduces to
\begin{eqnarray}\label{varrho3}
\frac{1}{\varepsilon k_1}(\varrho_t+A\varrho_x)+(\varrho u)_x=0.
\end{eqnarray}
To relate $\varrho$ to a physically meaningful variable, we set
\begin{eqnarray}\label{rho}
\rho^{\frac{1}{b-1}}=\varrho, \quad b\neq1.
\end{eqnarray}
Substituting (\ref{rho}) into (\ref{varrho3}), we get
\begin{eqnarray}\label{varrhotf}
{\rho _t}+A \rho_x+\varepsilon k_1\big((b-1)\rho u_x+u\rho_x\big)=0,
\end{eqnarray}
which, after the scaling (\ref{scaling1}) and a Galilean transformation, will yield
the second component of the R2$b$-family system (\ref{systemsus}).

On the other hand, it is inferred from (\ref{varrho}) retaining terms to order $O(\varepsilon^2,\varepsilon\mu)$ that
\begin{eqnarray*}\label{varrho4}
\varrho^2=1+\varepsilon (2k_1) \eta+\varepsilon^2(k_1^2 +2k_2) \eta^2+O(\varepsilon^3).
\end{eqnarray*}
Using the above equality, we express $\eta$ in terms of $\varrho$ to find
\begin{equation}\label{etatx1}
\left\{
\begin{aligned}
&\eta_t=\frac{1}{\varepsilon (2k_1)}(\varrho^2)_t-\varepsilon\frac{k_1^2 +2k_2}{2k_1}(\eta^2)_t+O(\varepsilon^2),\\
&\eta_x=\frac{1}{\varepsilon (2k_1)}(\varrho^2)_x-\varepsilon\frac{k_1^2 +2k_2}{2k_1}(\eta^2)_x+O(\varepsilon^2).
\end{aligned}
\right.
\end{equation}
Now we reexpress the $\varrho^2$-term in (\ref{etatx1}) in terms of $\rho^2$ by using
\begin{eqnarray}\label{rho1}
\varrho=1+\varepsilon\varrho_0,
\end{eqnarray}
where $\varrho_0$ is of order $O(1,\varepsilon,\mu)$. Retaining terms to order $O(\varepsilon^2)$, we apply a Taylor expansion to write
\begin{eqnarray*}\label{rho2}
\varrho^{2(b-1)}=(1+\varepsilon\varrho_0)^{2(b-1)}=1+2(b-1)\varepsilon\varrho_0+(b-1)(2b-3)\varepsilon^2\varrho^2_0+O(\varepsilon^3)=\rho^2,
\end{eqnarray*}
and as such
\begin{eqnarray}\label{rho3}
2\varepsilon\varrho_0=\frac{1}{b-1}\rho^2-(2b-3)\varepsilon^2\varrho^2_0-\frac{1}{b-1}+O(\varepsilon^3).
\end{eqnarray}
It follows the definition (\ref{rho1}) of $\varrho_0$ and (\ref{rho3}) that
\begin{eqnarray}\label{rho4}
\varrho^2=\frac{1}{b-1}\rho^2-(2b-4)\varepsilon^2\varrho^2_0+\frac{b-2}{b-1}+O(\varepsilon^3).
\end{eqnarray}
In addition, combining (\ref{varrho}) with (\ref{rho1}), it deduce from (\ref{rho4}) that $\varrho_0=k_1\eta +\varepsilon k_2\eta^2+O(\varepsilon^2)$, and then from the relations (\ref{relationueta}) that
\begin{eqnarray}\label{rho5}
\varrho^2&=&\frac{1}{b-1}\rho^2-(2b-4)\varepsilon^2k^2_1\eta^2+\frac{b-2}{b-1}+O(\varepsilon^3)\nonumber \\
&=&\frac{1}{b-1}\rho^2-\varepsilon^2
\frac{(2b-4)k^2_1}{(c-A)^2}u^2+\frac{b-2}{b-1}+O(\varepsilon^3).
\end{eqnarray}
The equation (\ref{rho5}) allows us to replace $\varrho^2$ appearing in (\ref{etatx1}) with the variable $\rho^2$. We have
\begin{equation}\label{etatx2}
\left\{
\begin{aligned}
&\eta_t=\frac{1}{\varepsilon (2k_1)}\big(\frac{1}{b-1}\rho^2-\varepsilon^2
\frac{(2b-4)k^2_1}{(c-A)^2}u^2\big)_t-\varepsilon\frac{k_1^2 +2k_2}{2k_1}(\eta^2)_t+O(\varepsilon^2),\\
&\eta_x=\frac{1}{\varepsilon (2k_1)}\big(\frac{1}{b-1}\rho^2-\varepsilon^2
\frac{(2b-4)k^2_1}{(c-A)^2}u^2\big)_x-\varepsilon\frac{k_1^2 +2k_2}{2k_1}(\eta^2)_x+O(\varepsilon^2).
\end{aligned}
\right.
\end{equation}
Then forming the linear combination $\eta_x+2\Omega\eta_t$ using the above expressions (\ref{etatx2}), we obtain
\begin{eqnarray}\label{Omegaeta}
&&\eta_x+2\Omega\eta_t\nonumber \\
&=&\frac{1}{(\varepsilon k_1)(b-1)}\rho(\rho_x+2\Omega\rho_t)-\varepsilon\frac{(2b-4)k_1}{(c-A)^2}u(u_x+2\Omega u_t)-\varepsilon\frac{k_1^2 +2k_2}{k_1}\eta(\eta_x+2\Omega \eta_t)+O(\varepsilon^2)\nonumber \\
&=&\frac{1}{(\varepsilon k_1)(b-1)}\rho(\rho_x+2\Omega\rho_t)-\varepsilon(2b-4)k_1\frac{1-2\Omega c}{(c-A)^2}uu_x-\varepsilon\frac{k_1^2 +2k_2}{k_1}\frac{c}{c-A}uu_x
+O(\varepsilon^2,\varepsilon\mu)\nonumber \\
&=&\frac{1}{(\varepsilon k_1)(b-1)}\rho(\rho_x+2\Omega\rho_t)-\varepsilon \big((2b-4)k_1+\frac{k_1^2 +2k_2}{k_1}
\big)\frac{c}{c-A}uu_x+O(\varepsilon^2,\varepsilon\mu),
\end{eqnarray}
where we have used $u_t=-cu_x+O(\varepsilon,\mu)$, ${\eta _x}+2\Omega\eta_t=-u_t=cu_x+O(\varepsilon,\mu)$ and the linear relations (\ref{relationueta}) in the above second equality.
We also used the dispersion relation $c^2-(A-2\Omega)c-1=0$ to obtain $\frac{1-2\Omega c}{(c-A)^2}=\frac{c}{c-A}$ in the above third equality.

Then in view of the the equation (\ref{varrhotf}) satisfied by $\rho$, substituting (\ref{Omegaeta}) into the second equation of (\ref{jianhuarGNwithshear}), we find
\begin{eqnarray}\label{u}
&&(u-\frac{\mu }{3}u_{xx})_t+ \varepsilon u{u_x}-\mu\frac{A}{3}u_{xxx}+\frac{1-2\Omega A}{(\varepsilon k_1)(b-1)}\rho\rho_x\nonumber \\
&-&\frac{2\Omega}{b-1}\rho\big((b-1)\rho u_x+u\rho_x\big)-\varepsilon \big((2b-4)k_1+\frac{k_1^2 +2k_2}{k_1}
\big)\frac{c}{c-A}uu_x=O(\varepsilon\mu,\mu^2).
\end{eqnarray}
Denoting $m=u-\frac{\mu }{3}u_{xx}$. Note that to the leading order we may break $uu_x$ up as
\begin{eqnarray*}\label{m}
uu_x=
\frac{\sigma}{1+b}(bu_xm+um_x)+(1-\sigma)uu_x+O(\mu),
\end{eqnarray*}
where $\sigma$ is an arbitrary real constant that originates from the balance between nonlinear steepening and stretching (see Introduction for its physical background).
Then the equation (\ref{u}) can be written at the order $O(\varepsilon,\mu)$ as
\begin{eqnarray}\label{u1}
&&m_t+Am_x-Au_x+\varepsilon\Big(1-\big((2b-3)k_1+\frac{2k_2}{k_1}
\big)\frac{c}{c-A}\Big)\frac{\sigma}{1+b}(bu_xm+um_x)\nonumber \\
&+&\varepsilon\Big(1-\big((2b-3)k_1+\frac{2k_2}{k_1}
\big)\frac{c}{c-A}\Big)(1-\sigma)uu_x+
\frac{1-2\Omega A}{(\varepsilon k_1)(b-1)}\rho\rho_x\nonumber \\
&-&\frac{2\Omega}{b-1}\rho\big((b-1)\rho u_x+u\rho_x\big)=0.
\end{eqnarray}
Having imposed only one constraint on the constants $k_1$ and $k_2$ in (\ref{relation1}), we now also choose $k_1$ and $k_2$ such that
\begin{eqnarray}\label{relation2}
1-\big((2b-3)k_1+\frac{2k_2}{k_1}
\big)\frac{c}{c-A}=k_1(b+1),
\end{eqnarray}
in which case the equation (\ref{u1}) becomes
\begin{eqnarray}\label{u2}
&&m_t+Am_x-Au_x+\varepsilon k_1\sigma(bu_xm+um_x)+\varepsilon k_1(b+1)(1-\sigma)uu_x\nonumber \\
&+&\frac{1-2\Omega A}{(\varepsilon k_1)(b-1)}\rho\rho_x-\frac{2\Omega}{b-1}\rho\big((b-1)\rho u_x+u\rho_x\big)=0.
\end{eqnarray}

Applying the following scaling
\begin{eqnarray}\label{scaling1}
x\rightarrow \sqrt{\frac{\mu}{3}}x,\quad t\rightarrow \sqrt{\frac{\mu}{3}}t,\quad u\rightarrow \frac{1}{k_1 \varepsilon}u, \quad \rho\rightarrow \sqrt{b-1}\rho,
\end{eqnarray}
to (\ref{varrhotf}) and (\ref{u2}), then we arrive at (After this scaling (\ref{scaling1}), we denote the new variables again by $u$, $\rho$, $m$ for simplicity; note that now $m = u - u_{xx}$, the factor $\mu/3$ originally present in $m = u - \frac{\mu}{3}u_{xx}$ having been absorbed into the scaling of $x$ and $t$.)
\begin{equation}\label{system}
\left\{
 \begin{aligned}
&m_t+Am_x-Au_x+\sigma(bu_xm+um_x)+(b+1)(1-\sigma)uu_x\\
&\quad\quad\quad\quad\quad\quad\quad\quad\ +(1-2\Omega A)\rho\rho_x-2\Omega\rho\big((b-1)\rho u_x+u\rho_x\big)=0,\\
&\rho_t+A\rho_x+u\rho_x+(b-1)\rho u_x=0,
\end{aligned}
\right.
\end{equation}
with the constants $k_1$, $k_2$ determined by the constraints (\ref{relation1}) and (\ref{relation2})
and $c$ satisfying the dispersion relation, namely
\begin{eqnarray}\label{k12}
k_1 = \frac{3c^2 - 3cA + A^2}{(c-A)\bigl(3bc - (b+1)A\bigr)},\
k_2 = k_1^2 - k_1\left(1 + \frac{A}{2(c-A)}\right),\  c^2-(A-2\Omega)c-1=0.
\end{eqnarray}
Under the Galilean transformation $x\rightarrow x-At$
, $t\rightarrow t$, the system (\ref{system}) becomes
\begin{equation*}
\left\{
 \begin{aligned}
&m_t-Au_x+\sigma(bu_xm+um_x)+(b+1)(1-\sigma)uu_x\\
&\quad\quad\quad\quad\quad\quad\quad\quad\quad\ +(1-2\Omega A)\rho\rho_x-2\Omega\rho\big((b-1)\rho u_x+u\rho_x\big)=0,\\
&\rho_t+u\rho_x+(b-1)\rho u_x=0.
\end{aligned}
\right.
\end{equation*}

We have now completed the derivation of the R2$b$-family system (\ref{systemsus}). In the following we comment on the physical meaning of the parameters and on some special reductions that follow from the algebraic relations (\ref{k12}).

\begin{remark2}
The following points are worth noting concerning the parameters in system (\ref{system}).
\begin{itemize}
\item The parameter $\sigma$ originates from the hyperelastic rod model of \cite{D,DH} and the generalized 2CH system \cite{CLiu}, where it controls the competition between the quadratic nonlinearity $uu_x$ and the higher-order dispersive terms $u_xu_{xx}$ and $uu_{xxx}$.  In our system \eqref{system}, $\sigma$ continues this role, appearing in the combination $\sigma(b u_x m + u m_x) + (b+1)(1-\sigma)uu_x$.
    \item The parameter $b$ plays the same role as in the single-component $b$-family equation \cite{HS,HS1}, where it governs the balance between convection and stretching.  In our system \eqref{system}, $b$ appears in the stretching term $(b-1)\rho u_x$ (second equation), in $\sigma(b u_x m + u m_x)$ and $(b+1)(1-\sigma)uu_x$ (first equation), and in the coupling term $2\Omega\rho\big((b-1)\rho u_x + u\rho_x\big)$.
        The condition $b\neq1$ is required by the definition $\rho^{\frac{1}{b-1}}=\varrho$ in (\ref{rho}).
        \item The inequality $1-2\Omega A>0$ is physically justified for realistic equatorial parameters ($\Omega\sim10^{-5}\,\text{rad}/\text{s}$, $A\sim10^{-2}$, water depth $h<300\,\text{m}$), as discussed in \cite{FGL,GQ}.
\end{itemize}
\end{remark2}

\begin{remark2}
If $k_2=0$, then $\varrho = 1 + k_1\varepsilon\eta$ and, using $\varrho = \rho^{1/(b-1)}$, we obtain
\[
\eta = \frac{\rho^{1/(b-1)} - 1}{k_1\varepsilon},
\]
which avoids solving a quadratic equation to recover $\eta$ from $\rho$ and may be advantageous in numerical simulations or analytical manipulations.
From the relations (\ref{k12}), the condition $k_2=0$ imposes a relation between $b$ and the physical parameters $A$ and $\Omega$
\[
b = \frac{6c^2 - 4cA + A^2}{(2c - A)(3c - A)}.
\]
Thus, for given $A$ and $\Omega$, the value of $b$ that makes $k_2$ vanish is uniquely determined.
\end{remark2}

\begin{remark2}
As discussed in Introduction, the R2$b$-family system (\ref{systemsus}) contains several known models as special cases:
\begin{itemize}
\item Rotation-two-component CH system (R2CH) \cite{FGL} when $b=2$;
\item Constant vorticity two-component $b$-family system \cite{EHKL} when $\Omega=0$, $\sigma=1$;
\item Generalized two-component CH system \cite{CLiu} when $b=2$, $\Omega=0$;
\item Single component $b$-family equation when $\rho\equiv0$ and $\sigma=1$ (with $A=0$ for the standard form).
\end{itemize}
\end{remark2}

\section{Wave-breaking}\label{sec3}
\newtheorem{theorem3}{Theorem}[section]
\newtheorem{lemma3}{Lemma}[section]
\newtheorem {remark3}{Remark}[section]
\newtheorem {definition3}{Definition}[section]
\newtheorem{corollary3}{Corollary}[section]
\newtheorem{example3}{Example}[section]
\par
In this section, we shall mainly consider the wave-breaking phenomena for the R2$b$-family system (\ref{systemsus}) with $b=3,A=0,\sigma=1,$ $i.e.,$ R2DP system (\ref{systems2}).
Note that one may apply Kato's semigroup theory \cite{Kato} as in \cite{FGL,YY} to obtain the following local well-posedness result for the Cauchy problem of the R2$b$-family system (\ref{systemsus}). Thus we omit the details here.
\begin{theorem3}\label{Th3.1}
Suppose that $(u_0,\rho_0-1)$ in $H^s\times H^{s-1}$ with $s>\frac{3}{2}$. Then there exist
 a maximal time $T=T(\|(u_0,\rho_0-1)\|_{H^s\times H^{s-1}})>0$ and a unique solution $(u,\rho-1) \in C([0,T);H^s\times H^{s-1})\cap C^1([0,T);H^{s-1}\times H^{s-2})$
 of the system (\ref{systemsus}) with $(u(0),\rho(0))=(u_0,\rho_0)$. Moreover, the solution $(u,\rho-1)$ depends continuously on the initial data $(u_0,\rho_0-1)$ and the existence time $T>$ can be chosen to be independent of $s$.
\end{theorem3}

In concern of the wave-breaking phenomena, we consider the initial value problem
\begin{equation}\label{flow1}
\begin{aligned}
\left\{ {\begin{array}{*{20}{l}}
q_{t}(t,x) =  u(t,q(t,x)),   &t \in [0,T),\\
q(0,x) = x,  & x \in \R,
\end{array}} \right.
\end{aligned}
\end{equation}
where $u\in C^1([0,T);H^{s-1})$ is the first component of the solution $(u,\rho)$ to the system (\ref{systemsus}) with initial data $(u_0,\rho_0)\in H^s\times H^{s-1}$ with $s>\frac{3}{2}$, and
$T>0$ is the maximal existence time. A direct calculation yields
\begin{eqnarray*}
q_{tx}(t,x)=u_{x}(t,q(t,x))q_{x}(t,x).
\end{eqnarray*}
Thus for $t>0,\ x\in \R$, we get
\begin{eqnarray}\label{q1}
q_{x}(t,x)=\exp\big(\int_0^tu_x(\tau,q(\tau,x))d\tau\big)>0,
\end{eqnarray}
which implies that the map $q(t,\cdot):\R \rightarrow \R$ is an increasing diffeomorphism of $\R$ for each $t\in [0,T).$ Hence for any function $v(t,\cdot)\in L^{\infty}(\R),t\in [0,T)$,
we have
\begin{eqnarray}\label{q2}
\|v(t,\cdot)\|_{L^{\infty}(\R)}=\|v(t,q(t,\cdot))\|_{L^{\infty}(\R)}, \quad t\in [0,T).
\end{eqnarray}

\subsection{Blow-up scenario}\label{Sec3.1}
In this subsection, we shall present the precise blow-up scenario for solutions of the R2DP system (\ref{systems2}).
Note that the second equation of the R2DP system (\ref{systems2}) is the same as the one of the 2DP system. As shown in \cite{YY},
The following lemma also holds true for (\ref{systems2}). We recall it for completeness.
\begin{lemma3}\label{Lem3.1}
(\cite{YY}) Let
$(u_0,\rho_0-1)\in H^s\times H^{s-1}$ with $s>\frac{3}{2}$ and $T>0$ be the maximal existence time of the corresponding solution $(u,\rho)$
 to (\ref{systems2}) (guaranteed by Theorem \ref{Th3.1}). Then
 \begin{eqnarray}\label{q3}
 \rho(t,q(t,x))\,q_x^2(t,x)=\rho_0(x),\qquad \forall (t,x)\in[0,T)\times\mathbb{R}.
 \end{eqnarray}
Moreover, if there exists $M>0$ such that $u_x(t,x)\ge -M$ for all $(t,x)\in[0,T)\times\mathbb{R}$, then
\[
\|\rho(t,\cdot)\|_{L^\infty}\le e^{2Mt}\|\rho_0\|_{L^\infty} \quad \mbox{and} \quad
\|\rho(t,\cdot)\|_{L^2}\le e^{2Mt}\|\rho_0\|_{H^{s-1}},\quad \forall t\in[0,T).
\]
\end{lemma3}
By using the local well-posedness result and energy estimates, we establish the precise blow-up scenario for sufficiently regular solutions to (\ref{systems2}) as fololows.
\begin{theorem3}\label{Th3.2}
Let $(u_0,\rho_0-1)$ in $H^s\times H^{s-1}$ with $s>\frac{5}{2}$, and let
$T>0$ be the maximal existence time of the solution $(u,\rho)$ to the R2DP system (\ref{systems2}) with initial data $(u_0,\rho_0)$.
Then the corresponding solution blows up in finite time if and only if
 \begin{eqnarray*}
 \lim_{ t\rightarrow T^-}\inf_{x\in\R}\{u_x(t,x)\}=-\infty\ \  \mbox{or}\ \ \limsup_{ t\rightarrow T^-}\{\|u(t,\cdot)\|_{L^\infty}\}=+\infty,\ \  \mbox{or}\ \ \limsup_{ t\rightarrow T^-}\{\|\rho_{x}(t,\cdot)\|_{L^\infty}\}=+\infty.
 \end{eqnarray*}
\end{theorem3}
\begin{proof}
Denote $\bar{\rho}:=\rho-1$, then we can rewrite the system (\ref{systems2}) as
\begin{equation}\label{systems4}
\left\{
 \begin{aligned}
&m_t+um_x+3u_xm+\bar{\rho}\bar{\rho}_x+\bar{\rho}_x-2\Omega\big(2
\bar{\rho}^2u_x+4\bar{\rho} u_x+u\bar{\rho}\bar{\rho}_x+2u_x+u\bar{\rho}_x
\big)=0,\\
&\bar{\rho}_t+u\bar{\rho}_x+2\bar{\rho} u_x+2u_x=0.
\end{aligned}
\right.
\end{equation}
Multiplying the second equation in (\ref{systems4}) by $\bar{\rho}$ and integrating by parts, we have
 \begin{eqnarray}\label{3.6}
 \frac{1}{2}\frac{d}{dt} \int_{\R} \bar{\rho}^2dx=-\frac{3}{2} \int_{\R} u_x\bar{\rho}^2dx-2 \int_{\R} \bar{\rho} u_xdx.
\end{eqnarray}
Differentiating the second equation in (\ref{systems4}) with respect to $x$,
multiplying the result equation by $\bar{\rho}_x$ and integrating by parts, we get
 \begin{eqnarray}\label{3.7}
 \frac{1}{2}\frac{d}{dt} \int_{\R} \bar{\rho}_x^2dx=-\frac{5}{2} \int_{\R} u_x\bar{\rho}_x^2dx+\int_{\R} \bar{\rho}^2 u_{xxx}dx
 -2 \int_{\R} \bar{\rho}_x u_{xx}dx.
\end{eqnarray}
Differentiating the second equation in (\ref{systems4}) with respect to $x$ twice,
multiplying the result equation by $\bar{\rho}_{xx}$ and integrating by parts, we obtain
 \begin{eqnarray}\label{3.8}
 \frac{1}{2}\frac{d}{dt} \int_{\R} \bar{\rho}_{xx}^2dx=-\frac{7}{2} \int_{\R} u_x\bar{\rho}_{xx}^2dx+\frac{5}{2} \int_{\R} u_{xxx}\bar{\rho}_x^2dx
 -2 \int_{\R} \bar{\rho}\bar{\rho}_{xx} u_{xxx}dx-2 \int_{\R}\bar{\rho}_{xx} u_{xxx}dx.
\end{eqnarray}
On the other hand, multiplying the first equation in (\ref{systems4}) by $m=u-u_{xx}$ and integrating by parts, we have
 \begin{eqnarray}\label{3.9}
 \frac{1}{2}\frac{d}{dt} \int_{\R} m^2dx&=&-\frac{5}{2}\int_{\R} u_x m^2dx+\frac{1}{2} \int_{\R} \bar{\rho}^2 u_xdx-\frac{1}{2} \int_{\R} \bar{\rho}^2 u_{xxx}dx
 +\int_{\R}\bar{\rho} u_xdx+\int_{\R}\bar{\rho}_x u_{xx}dx\nonumber\\
 &&-2\Omega\int_{\R}\bar{\rho}\bar{\rho}_x\big(u^2-2u^2_x +uu_{xx}\big)dx-2\Omega\int_{\R}\bar{\rho}_x\big(u^2-2u^2_x+uu_{xx}\big)dx.\nonumber\\
\end{eqnarray}
Differentiating the first equation in (\ref{systems4}) with respect to $x$,
multiplying the result equation by $m_x=u_x-u_{xxx}$ and integrating by parts, we get
\begin{eqnarray}\label{3.10}
 \frac{1}{2}\frac{d}{dt} \int_{\R} m_x^2dx&=&-\frac{7}{2}\int_{\R} u_x m_x^2dx+\frac{3}{2}\int_{\R} u_x m^2dx
 +\int_{\R}u_{xxx}(\bar{\rho}_x^2+\bar{\rho}\bar{\rho}_{xx}-\frac{1}{2}\bar{\rho}^2)dx+\int_{\R}\bar{\rho}_x u_{xx}dx\nonumber\\
 &&+\int_{\R}\bar{\rho}_{xx} u_{xxx}dx+2\Omega\int_{\R}\bar{\rho}\bar{\rho}_x\big(2u^2_x-5 u_x u_{xx}-2u^2_{xx}-uu_{xx}\big)dx\nonumber\\
 &&-2\Omega\int_{\R}\bar{\rho}^2_xuu_{xxx}dx+2\Omega\int_{\R}\bar{\rho}_x\big(2u^2_x-u u_{xx}+2u^2_{xx}-u_xu_{xxx}\big)dx\nonumber\\
 &&-2\Omega\int_{\R}u\bar{\rho}_{xx}u_{xxx}dx-2\Omega\int_{\R}u\bar{\rho}\bar{\rho}_{xx}u_{xxx}dx.
\end{eqnarray}
Combining (\ref{3.6})-(\ref{3.8}) with (\ref{3.9})-(\ref{3.10}), we find
\begin{eqnarray}\label{3.11}
&&\frac{d}{dt} \int_{\R}\big(\frac{1}{2}\bar{\rho}^2+\bar{\rho}_x^2+\frac{1}{2}\bar{\rho}_{xx}^2
+m^2+m_x^2\big)dx\nonumber\\
&=&- \int_{\R}u_x\big(\frac{1}{2}\bar{\rho}^2+5\bar{\rho}_x^2+\frac{7}{2}\bar{\rho}_{xx}^2+2m^2+7m_x^2\big)dx+\frac{9}{2}\int_{\R}u_{xxx}\bar{\rho}_x^2dx\nonumber\\
&&+4\Omega\int_{\R}\bar{\rho}\bar{\rho}_{x}\big(-u^2+4u_x^2-2uu_{xx}-5u_xu_{xx}-2u^2_{xx}\big)dx-4\Omega\int_{\R}u u_{xxx}\bar{\rho}_x^2dx\nonumber\\
&&+4\Omega\int_{\R}\bar{\rho}_{x}\big(-u^2+4u_x^2-2uu_{xx}+2u^2_{xx}-u_xu_{xxx}\big)dx\nonumber\\
&&-4\Omega\int_{\R}(1+\bar{\rho})\bar{\rho}_{xx}uu_{xxx}dx.
\end{eqnarray}
Assume that $T<+\infty$ and there exists $M>0$ such that
 \begin{eqnarray}\label{3.12}
 u_x(t,x)\geq -M, \ \forall(t,x)\in[0,T)\times \R, \ \mbox{and} \
 \|u(t,\cdot)\|_{L^\infty}, \|\rho_{x}(t,\cdot)\|_{L^\infty}\leq M, \ \forall t\in[0,T).
 \end{eqnarray}
Hence, in view of Lemma \ref{Lem3.1},  we also have
\begin{eqnarray}\label{3.13}
 \|\rho(t,\cdot)\|_{L^\infty}\le e^{2MT}\|\rho_0\|_{L^\infty},\qquad \forall t\in[0,T).
 \end{eqnarray}
Consequently, from (\ref{3.11})-(\ref{3.13}),
there exists a constant $C(M,\Omega,T,\|\rho_0\|_{L^\infty})>0$, depending only on $M$, $\Omega$, $T$, and $\|\rho_0\|_{L^\infty}$, such that
\begin{eqnarray*}
&&\frac{d}{dt} \int_{\R}\big(\frac{1}{2}\bar{\rho}^2+\bar{\rho}_x^2+\frac{1}{2}\bar{\rho}_{xx}^2
+m^2+m_x^2\big)dx\nonumber\\
&\leq& C(M,\Omega,T,\|\rho_0\|_{L^\infty}) \big(\frac{1}{2}\bar{\rho}^2+\bar{\rho}_x^2+\frac{1}{2}\bar{\rho}_{xx}^2
+m^2+m_x^2\big)dx.
\end{eqnarray*}
Applying Gronwall's inequality to the above inequality yields for every $t\in [0,T)$
\begin{eqnarray*}
\|u\|^2_{H^3}+\|\bar{\rho}\|^2_{H^2}\leq \|m\|^2_{H^1}+\|\bar{\rho}\|^2_{H^2}
 \leq 2 (\|m(0,\cdot)\|^2_{H^1}+\|\bar{\rho}(0)\|^2_{H^2})e^{C(M,\Omega,T,\|\rho_0\|_{L^\infty})t}, \quad \forall t\in[0,T),
\end{eqnarray*}
which contradicts the assumption the maximal existence time $T<+\infty$.

Conversely, by Sobolev embedding theorem $H^s(\R)\hookrightarrow L^\infty(\R)$, $s>\frac{1}{2}$, we see that if
\begin{eqnarray*}
 \lim_{ t\rightarrow T^-}\inf_{x\in\R}\{u_x(t,x)\}=-\infty\ \  \mbox{or}\ \ \limsup_{ t\rightarrow T^-}\{\|u(t,\cdot)\|_{L^\infty}\}=+\infty,\ \  \mbox{or}\ \ \limsup_{ t\rightarrow T^-}\{\|\rho_{x}(t,\cdot)\|_{L^\infty}\}=+\infty,
 \end{eqnarray*}
then the solution blows up in finite time. This completes the proof of Theorem \ref{Th3.2}.
\end{proof}

Indeed, for solutions with lower regularity, a simplified energy estimate suffices.
Combining the lower-order identities in the proof of  Theorem \ref{Th3.2}, we have the following blow-up scenario.
\begin{theorem3}\label{Th3.3}
Let $(u_0,\rho_0-1)$ in $H^2\times H^{1}$, and let
$T>0$ be the maximal existence time of the solution $(u,\rho)$ to the R2DP system (\ref{systems2}) with initial data $(u_0,\rho_0)$.
Then the corresponding solution blows up in finite time if and only if
 \begin{eqnarray*}
 \lim_{ t\rightarrow T^-}\inf_{x\in\R}\{u_x(t,x)\}=-\infty\ \  \mbox{or}\ \ \limsup_{ t\rightarrow T^-}\{\|u(t,\cdot)\|_{L^\infty}\}=+\infty.
 \end{eqnarray*}
\end{theorem3}
\begin{proof}
Form (\ref{3.6})-(\ref{3.7}) and (\ref{3.9}), we obtain
\begin{eqnarray}\label{3.14}
&&\frac{d}{dt} \int_{\R}\big(\frac{1}{2}\bar{\rho}^2+\frac{1}{2}\bar{\rho}_x^2+m^2\big)dx\nonumber\\
&=&- \int_{\R}u_x\big(\frac{1}{2}\bar{\rho}^2+\frac{5}{2}\bar{\rho}_x^2+5m^2\big)dx+4\Omega\int_{\R}\bar{\rho}^2\big(uu_x-2u_xu_{xx}\big)dx\nonumber\\
&&+8\Omega\int_{\R}\bar{\rho}(uu_x-2u_xu_{xx})dx-4\Omega\int_{\R}(\bar{\rho}+1)u \bar{\rho}_x u_{xx}dx.
\end{eqnarray}
Suppose, for contradiction, that $T<+\infty$ and there exists $M>0$ such that
 \begin{eqnarray}\label{3.15}
 u_x(t,x)\geq -M, \ \forall(t,x)\in[0,T)\times \R, \ \mbox{and} \
 \|u(t,\cdot)\|_{L^\infty}\leq M, \ \forall t\in[0,T).
 \end{eqnarray}
Thus, by Lemma \ref{Lem3.1}, we also have
\begin{eqnarray}\label{3.16}
 \|\rho(t,\cdot)\|_{L^\infty}\le e^{2MT}\|\rho_0\|_{L^\infty},\qquad \forall t\in[0,T).
 \end{eqnarray}
Therefore, using (\ref{3.14})-(\ref{3.16}),
there exists a constant $C(M,\Omega,T,\|\rho_0\|_{L^\infty})>0$, depending only on $M$, $\Omega$, $T$, and $\|\rho_0\|_{L^\infty}$, such that
\begin{eqnarray}\label{3.17}
\frac{d}{dt} \int_{\R}\big(\frac{1}{2}\bar{\rho}^2+\frac{1}{2}\bar{\rho}_x^2+m^2\big)dx
\leq C(M,\Omega,T,\|\rho_0\|_{L^\infty})\int_{\R}\big(\frac{1}{2}\bar{\rho}^2+\frac{1}{2}\bar{\rho}_x^2+m^2\big)dx.
\end{eqnarray}
Applying Gronwall's inequality to (\ref{3.17}) yields for every $t\in [0,T)$
\begin{eqnarray*}
\|u\|^2_{H^2}+\|\bar{\rho}\|^2_{H^1}\leq \|m\|^2_{L^2}+\|\bar{\rho}\|^2_{H^1}
 \leq 2(\|m(0,\cdot)\|^2_{L^2}+\|\bar{\rho}(0)\|^2_{H^1})e^{C(M,\Omega,T,\|\rho_0\|_{L^\infty})
t}, \quad \forall t\in[0,T),
\end{eqnarray*}
which contradicts the assumption that the maximal existence time $T<+\infty$.

Conversely, if
\begin{eqnarray*}
 \lim_{ t\rightarrow T^-}\inf_{x\in\R}\{u_x(t,x)\}=-\infty\ \  \mbox{or}\ \ \limsup_{ t\rightarrow T^-}\{\|u(t,\cdot)\|_{L^\infty}\}=+\infty,
 \end{eqnarray*}
then by Sobolev embedding theorem $H^s(\R)\hookrightarrow L^\infty(\R)$, $s>\frac{1}{2}$, the solution cannot exist beyond $T$, hence blows up in finite time.
This completes the proof of Theorem \ref{Th3.3}.
\end{proof}

\subsection{Blow-up}\label{Sec3.2}
In this subsection, we shall address the question of the formation of singularity for solutions of the R2DP system (\ref{systems2}).
As discussed in Introduction, for the R2DP system (\ref{systems2}) the absence of a conserved energy (as in the 2DP system \cite{YY}) and the presence of the cubic coupling term $\rho(2\rho u_x+u\rho_x)$ (as in the R2CH system \cite{FGL}) render the classical blow-up criteria inapplicable.  To overcome this double difficulty, we introduce an auxiliary variable that absorbs the most singular part of the coupling term.  Following the idea in \cite{FGL}, we define
 \begin{eqnarray}\label{omega}
 \omega= u +\Omega\,p\ast\rho^2,\qquad p(x)=\frac12 e^{-|x|}.
 \end{eqnarray}
Then the evolution equation for $\omega$ is given in the following lemma.
\begin{lemma3}\label{Lem3.2}
For the R2DP system (\ref{systems2}), the auxiliary variable $\omega$ given by (\ref{omega}) satisfies
 \begin{eqnarray}\label{K}
\omega_t+u\omega_x
=\Omega \omega \partial_xp\ast (\rho^2)-\Omega^2(p\ast\rho^2)(\partial_xp\ast\rho^2)-\partial_xp\ast (\frac{3}{2}u^2+\frac{1}{2}\rho^2),
 \end{eqnarray}
 and
  \begin{eqnarray}\label{Ktx}
\omega_{tx}+u\omega_{xx}
=-\big(\omega_x-\Omega\,\partial_xp\ast\rho^2\big)^2+(\frac{1}{2}-\Omega u)\rho^2+\Omega u p\ast\rho^2
+\frac{3}{2}u^2-p\ast(\frac{3}{2}u^2+\frac{1}{2}\rho^2),
 \end{eqnarray}
 where $p(x)=\frac{1}{2}e^{-|x|},x\in \R$.
\end{lemma3}

\begin{proof}
From the second equation of (\ref{systems2}) we have $\rho_t = -2\rho u_x - u\rho_x$, which implies
 \begin{eqnarray}\label{rhorhot}
\rho(2\rho u_x + u\rho_x) = -\rho\rho_t.
 \end{eqnarray}
Substituting (\ref{rhorhot}) into the first equation of (\ref{systems2}) gives
\begin{eqnarray}\label{uut}
u_t - u_{txx} +2\,\Omega\,\rho\rho_t =-4uu_x +3u_xu_{xx} + uu_{xxx} - \rho\rho_x.
 \end{eqnarray}
Using the identity $3u_xu_{xx} + uu_{xxx} = -(1-\partial_x^2)(uu_x) + uu_x$, we rewrite (\ref{uut}) as
\begin{eqnarray}\label{uut1}
\big((1-\partial_x^2)u+\Omega\rho^2\big)_t=-(1-\partial_x^2)(uu_x) -\frac{3}{2}(u^2)_x-\frac{1}{2}(\rho^2)_x.
 \end{eqnarray}
Applying the operator $(1-\partial_x^2)^{-1}$ to both sides of (\ref{uut1}), we obtain
\begin{eqnarray*}
&&\big(u +\Omega (1-\partial_x^2)^{-1} \rho^2\big)_t
=-(uu_x)-\partial_x(1-\partial_x^2)^{-1} (\frac{3}{2}u^2+\frac{1}{2}\rho^2)\nonumber \\
&=&-u\big(u + \Omega (1-\partial_x^2)^{-1} \rho^2\big)_x+\Omega u (1-\partial_x^2)^{-1} (\rho^2)_x-\partial_x(1-\partial_x^2)^{-1} (\frac{3}{2}u^2+\frac{1}{2}\rho^2).
 \end{eqnarray*}
In light of $\omega = u +\Omega\,p\ast\rho^2$, and noting that $(1-\partial^2_x)^{-1}f=p\ast f$ for all $f\in L^2(\R)$, we find
\begin{eqnarray}\label{uut2}
\omega_t+u\omega_x=\Omega u \partial_xp\ast (\rho^2)-\partial_xp\ast (\frac{3}{2}u^2+\frac{1}{2}\rho^2),
 \end{eqnarray}
which together with (\ref{omega}) proves (\ref{K}).

As for (\ref{Ktx}), differentiating (\ref{uut2}) with respect to $x$, and using $\partial_x^2p\ast f=p\ast f-f$, we have
\begin{eqnarray}\label{uut3}
\omega_{tx}+u\omega_{xx}&=&-u_x\omega_x+\Omega u_x\partial_xp\ast (\rho^2)+\Omega u\partial^2_xp\ast (\rho^2)-\partial^2_xp\ast (\frac{3}{2}u^2+\frac{1}{2}\rho^2)\nonumber\\
&=&-u^2_x+(\frac{1}{2}-\Omega u)\rho^2+\Omega u p\ast\rho^2
+\frac{3}{2}u^2-p\ast(\frac{3}{2}u^2+\frac{1}{2}\rho^2),
 \end{eqnarray}
which, combined with (\ref{omega}), gives (\ref{Ktx}).
\end{proof}

In the next lemma, we establish a priori estimate for the $L^\infty$-norm of the first component $u$ of the R2DP system (\ref{systems2}).
\begin{lemma3}\label{Lem3.3}
Let
$(u_0,\rho_0-1)\in H^s\times H^{s-1}$ with $s>\frac{3}{2}$ and $T>0$ be the maximal existence time of the corresponding solution $(u,\rho)$
 to (\ref{systems2}) (guaranteed by Theorem \ref{Th3.1}).
Suppose there exists a constant $M>0$ such that $\|\rho(t,\cdot)\|_{L^{2}}\le e^{2Mt}\|\rho_0\|_{H^{s-1}}$ for all $t\in[0,T)$.
Let $\omega$ be defined by (\ref{omega}). Then there exists a constant $\kappa>0$, depending only on $\|\omega_0\|_{L^\infty}$, $\|\rho_0\|_{H^{s-1}}$, $\Omega$, and $M$, such that for all
\begin{equation*}\label{time}
0\le t \le T^*:=\min\{T,\kappa\},
\end{equation*}
the following a priori estimate holds:
\begin{equation}\label{estimateu1}
\|u(t)\|_{L^\infty}
\le
\frac{\|\omega_0\|_{L^\infty} + A_0}
{1 - C_0\bigl(\|\omega_0\|_{L^\infty} + A_0\bigr)t}
- \frac{1}{6}\Omega R_\kappa^2,
\end{equation}
where
\[
R_\kappa := e^{2M\kappa}\|\rho_0\|_{H^{s-1}},\qquad
A_0 := \frac{2\Omega}{3}R_\kappa^2,\qquad
C_0 := \frac{3}{2} + \frac{1}{4}R_\kappa^2,
\]
and $\kappa$ is chosen such that
\begin{equation*}
\kappa < \frac{1}{C_0\bigl(\|\omega_0\|_{L^\infty} + A_0\bigr)}.
\end{equation*}
\end{lemma3}

\begin{proof}
Applying Theorem \ref{Th3.1} and a simple density argument, here we only need to prove the above lemma with $s \geq 3$.
Along the trajectory of $q(t,x)$ ($q(t,x)$ is defined in (\ref{flow1})), we deduce from equation (\ref{K}) that
 \begin{eqnarray}\label{Kflow}
\frac{d\omega}{dt}
=\Omega \omega \partial_xp\ast (\rho^2)-\Omega^2(p\ast\rho^2)(\partial_xp\ast\rho^2)-\partial_xp\ast (\frac{3}{2}u^2+\frac{1}{2}\rho^2),
 \end{eqnarray}
at $(t,q(t,x))$. Noting that
$\left|\partial_x p\ast (\frac{3}{2}u^2) \right| \le p \ast \frac{3}{2}u^2$, and applying Young's inequality, we obtain
 \begin{eqnarray}\label{estimate1}
\|\Omega \omega \partial_xp\ast (\rho^2)\|_{L^\infty}\leq \Omega \|\omega\|_{L^\infty}\| \partial_xp\|_{L^\infty}\|\rho^2\|_{L^1}\leq
 \frac{\Omega}{2} \big(e^{2Mt}\|\rho_0\|_{H^{s-1}}\big)^2 \|\omega\|_{L^\infty},
 \end{eqnarray}

 \begin{eqnarray}\label{estimate2}
\|\Omega^2(p\ast\rho^2)(\partial_xp\ast\rho^2)\|_{L^\infty}\leq \Omega^2 \|p\|_{L^\infty} \|\partial_xp\|_{L^\infty}\|\rho^2\|^2_{L^1}\leq
 \frac{\Omega^2}{4} \big(e^{2Mt}\|\rho_0\|_{H^{s-1}}\big)^4 ,
 \end{eqnarray}

 \begin{eqnarray}\label{estimate3}
\|\frac{1}{2}\partial_xp\ast (\rho^2)\|_{L^\infty}\leq \frac{1}{2} \|\partial_xp\|_{L^\infty}\|\rho^2\|_{L^1}\leq
 \frac{1}{4} \big(e^{2Mt}\|\rho_0\|_{H^{s-1}}\big)^2 ,
 \end{eqnarray}
and
\begin{eqnarray}\label{estimate4}
&&\|\frac{3}{2}\partial_xp\ast (u^2)\|_{L^\infty}\leq \|\frac{3}{2} p\ast (u^2)\|_{L^\infty}\leq
\frac{3}{2} \|p\|_{L^1}\|u^2\|_{L^\infty}\leq
 \frac{3}{2} \|u\|^2_{L^\infty}\nonumber\\
 &\leq& \frac{3}{2} \big(\|\omega\|_{L^\infty}+\|\Omega  p\ast (\rho^2)\|_{L^\infty}
\big)^2\leq  \frac{3}{2} \big(\|\omega\|_{L^\infty}+ \frac{\Omega}{2} (e^{2Mt}\|\rho_0\|_{H^{s-1}})^2
\big)^2.
 \end{eqnarray}
Combining (\ref{Kflow}) with these estimates (\ref{estimate1})-(\ref{estimate4}), we find
 \begin{eqnarray}\label{Kflow1}
\left| \frac{d\omega}{dt}\right|\leq \frac{3}{2}\big(\|\omega\|_{L^\infty}+\frac{2\Omega}{3}(e^{2Mt}\|\rho_0\|_{H^{s-1}})^2\big)^2+\frac{1}{4} \big(e^{2Mt}\|\rho_0\|_{H^{s-1}}\big)^2.
 \end{eqnarray}

Fix a positive number $\kappa$ (to be determined below) and restrict $t$ to
$0\le t\le T^*=\min\{T,\kappa\}$. Define
\begin{equation*}
R_\kappa := e^{2M\kappa}\|\rho_0\|_{H^{s-1}}.
\end{equation*}
Then on $[0,T^*]$, $R(t)=e^{2Mt}\|\rho_0\|_{H^{s-1}}\le R_\kappa$. Combining this with the differential inequality (\ref{Kflow1}), we have
 \begin{equation}\label{Kflow2}
\left| \frac{d\omega}{dt}\right|\le \frac{3}{2}\left(\|\omega\|_{L^\infty}+\frac{2\Omega}{3}R_\kappa^2\right)^2+\frac{1}{4}R_\kappa^2.
\end{equation}
Integrating (\ref{Kflow2}) along the characteristic $q(t,x)$ and taking the supremum over $x\in\mathbb{R}$ (since $q(t,\cdot)$ is a diffeomorphism) yields
\begin{equation}\label{Kflow3}
y(t)\le y(0)+\int_0^t \left[\frac{3}{2}(y(s)+A_0)^2+B_0\right]ds,
\end{equation}
where we set
\begin{equation*}
y(t):=\|\omega(t)\|_{L^\infty},\qquad
A_0:=\frac{2\Omega}{3}R_\kappa^2,\qquad
B_0:=\frac{1}{4}R_\kappa^2.
\end{equation*}
Define $z(t):=y(t)+A_0$. Then $z(0)=y(0)+A_0$, and (\ref{Kflow3}) becomes
\begin{equation*}
z(t)\le z(0)+\int_0^t \left[\frac{3}{2}z(s)^2+B_0\right]ds.
\end{equation*}
Without loss of generality, we assume $z(s)\ge 1$, then
\begin{equation*}
z(t)\le z(0)+C_0\int_0^t z(s)^2\,ds.
\end{equation*}
where $C_0:=\frac{3}{2}+B_0.$
Define $F(t):=z(0)+C_0\int_0^t z(s)^2\,ds$. Then $z(t)\le F(t)$, $F(0)=z(0)$, and
\[
F'(t)=C_0 z(t)^2 \le C_0 F(t)^2.
\]
Solving the above differential inequality gives
\[
F(t)\le \frac{z(0)}{1-C_0 z(0)t},\qquad \text{for } \ t<\frac{1}{C_0 z(0)}.
\]
Consequently,
\[
z(t)\le \frac{z(0)}{1-C_0 z(0)t},
\]
and hence
\begin{equation}\label{Kflow4}
y(t)\le \frac{y(0)+A_0}{1-C_0(y(0)+A_0)t}-A_0.
\end{equation}
To ensure that (\ref{Kflow4}) holds on the whole interval $[0,T^*]$, we require
\begin{equation}\label{Kflow5}
\kappa<\frac{1}{C_0(y(0)+A_0)}.
\end{equation}
Notice that $y(0)=\|\omega_0\|_{L^\infty}$, and both $A_0$ and $C_0$ depend on $R_\kappa=e^{2M\kappa}\|\rho_0\|_{H^{s-1}}$. Thus (\ref{Kflow5}) is an implicit condition on $\kappa$. However, by continuity, as $\kappa\to 0^+$, the right-hand side of (\ref{Kflow5}) tends to
\[
\frac{1}{\left(\frac{3}{2}+\frac{1}{4}\|\rho_0\|_{H^{s-1}}^2\right)\left(\|\omega_0\|_{L^\infty}+\frac{2\Omega}{3}\|\rho_0\|_{H^{s-1}}^2\right)}>0,
\]
so there exists at least one $\kappa>0$ satisfying (\ref{Kflow5}). Choosing such an $\kappa$, then for all $0\le t\le T^*$, the estimate (\ref{Kflow4}) is valid.

Finally, since $\omega = u + \Omega p*\rho^2$ and $\|p*\rho^2\|_{L^\infty}\le \frac{1}{2}R_\kappa^2$, we have
\[
\|u(t)\|_{L^\infty}\le \|\omega(t)\|_{L^\infty}+\frac{1}{2}\Omega R_\kappa^2.
\]
Substituting (\ref{Kflow4}) into the above inequality yields the desired estimate (\ref{estimateu1}). This completes the proof of Lemma \ref{Lem3.3}.
\end{proof}

By Lemma \ref{Lem3.1}, in order to study the fine structure of finite time singularities, we shall assume in the following
that there exists $M>0$ such that $\|\rho(t,\cdot)\|_{L^\infty},\ \|\rho(t,\cdot)\|_{L^2}\le e^{2Mt}\|\rho_0\|_{H^{s-1}}$ for all $t \in [0,T)$.

\begin{theorem3}\label{Th3.4}
 Let
$(u_0,\rho_0-1)\in H^s\times H^{s-1}$ with $s>\frac{3}{2}$ and $T>0$ be the maximal existence time of the corresponding solution $(u,\rho)$
 to (\ref{systems2}) (guaranteed by Theorem \ref{Th3.1}). Assume that there exists $M>0$ such that $\|\rho(t,\cdot)\|_{L^\infty},\ \|\rho(t,\cdot)\|_{L^2}\le e^{2Mt}\|\rho_0\|_{H^{s-1}}$ for all $t \in [0,T)$.
For some $\kappa_0\in(0,\kappa]$ ($\kappa$ is defined in Lemma \ref{Lem3.3}) and if there exists a point $x_0\in \R$ such that
 \begin{eqnarray}\label{x0condition}
 u_{0,x}(x_0)\leq -\sqrt{2\nu}-\frac{2\sqrt{2\nu}}{e^{\sqrt{2\nu}\kappa_0}-1}-\frac{3}{4}A_0,
 \end{eqnarray}
then corresponding solution $(u,\rho)$
 to (\ref{systems2}) blows up in finite time in the following sense, there is a $T$ with
 \begin{eqnarray*}
T\leq \frac{1}{\sqrt{2\nu}}\ln{\frac{u_{0,x}+\Omega\,\partial_{0,x}p\ast\rho^2(x_0)
-\sqrt{2\nu}}{u_{0,x}+\Omega\,\partial_{0,x}p\ast\rho^2(x_0)+\sqrt{2\nu}}}\leq \kappa_0,
 \end{eqnarray*}
 such that
  \begin{eqnarray*}
\liminf_{t\rightarrow T^{-}}\{\inf_{x\in\R}u_x(t,x)\}=-\infty,
 \end{eqnarray*}
 where $\nu:=\big(\frac{\Omega^2}{4} +\frac{1}{2}+\frac{3\Omega}{2} U(\kappa_0)\big)R_{\kappa_0}^2+\frac{3}{2}U^2(\kappa_0)$
 with $U(\kappa_0)=\frac{\|\omega_0\|_{L^\infty} + A_0}
{1 - C_0\bigl(\|\omega_0\|_{L^\infty} + A_0\bigr)\kappa_0}
- \frac{1}{6}\Omega R_{\kappa_0}^2$, $R_{\kappa_0}= e^{2M{\kappa_0}}\|\rho_0\|_{H^{s-1}}$
and $A_0, C_0$ are defined as in Lemma \ref{Lem3.3} with \(\kappa\) replaced by \(\kappa_0\).
\end{theorem3}
\begin{proof}
Applying Theorem \ref{Th3.1} and a simple density argument, we need only to show that Theorem \ref{Th3.4} holds with some $s\geq3$.
Along the trajectory of $q(t,x)$,
applying the inequality $(a-b)^2\ge \frac12 a^2-b^2$ and neglecting the non-positive terms $-p*(\frac32 u^2+\frac12\rho^2)\le 0$, we deduce from (\ref{Ktx}) that
\begin{eqnarray}\label{Ktx1}
\frac{d\omega_{x}(t,q(t,x))}{dt}
&\leq& -\frac{\omega^2_{x}(t,q(t,x))}{2}+f(t,q(t,x)),
 \end{eqnarray}
where $f:=\Omega^2(\partial_xp\ast\rho^2)^2+(\frac{1}{2}-\Omega u)\rho^2+\Omega u p\ast\rho^2
+\frac{3}{2}u^2$.
It then follows from the a priori estimate (\ref{estimateu1}), Young's inequality and the assumptions on $\rho$, we deduce for all $t\in [0,\kappa_0]\cap [0,T)$
\begin{eqnarray}\label{festimate}
f &\leq&  \Omega^2 (\|\partial_xp\|_{L^\infty}\|\rho^2\|_{L^1})^2+(\frac{1}{2}+\Omega\|u\|_{L^\infty})\|\rho\|^2_{L^\infty}+\Omega\|u\|_{L^\infty} \|p\|_{L^\infty}\|\rho^2\|_{L^1}+\frac{3}{2}\|u\|^2_{L^\infty}\nonumber \\
&\leq& \big(\frac{\Omega^2}{4}R_{\kappa_0}^2 +\frac{1}{2}+\frac{3\Omega}{2} U(\kappa_0)\big)R_{\kappa_0}^2+\frac{3}{2}U^2(\kappa_0):=\nu,
 \end{eqnarray}
where $U(\kappa_0)=\frac{\|\omega_0\|_{L^\infty} + A_0}
{1 - C_0\bigl(\|\omega_0\|_{L^\infty} + A_0\bigr)\kappa_0}
- \frac{1}{6}\Omega R_{\kappa_0}^2$, $R_{\kappa_0}= e^{2M{\kappa_0}}\|\rho_0\|_{H^{s-1}}$
and $A_0, C_0$ are defined as in Lemma \ref{Lem3.3} with \(\kappa\) replaced by \(\kappa_0\).
Hence, combining (\ref{Ktx1}) with (\ref{festimate}), and denoting $\bar{\omega}_{x}(t):=\omega_{x}(t,q(t,x_0))$,
we obtain for all $t\in [0,\kappa_0]\cap [0,T)$
\begin{eqnarray}\label{Ktx2}
\frac{d\bar{\omega}_{x}(t)}{dt}
&\leq& -\frac{\bar{\omega}_{x}(t)}{2}+\nu.
 \end{eqnarray}

 Now we claim that
 \begin{eqnarray}\label{Ktx3}
\bar{\omega}_{x}(t)=\omega_{x}(t,q(t,x_0))\leq \omega_{0,x}(x_0):=\bar{\omega}_x(0)<-\sqrt{2\nu},
 \end{eqnarray}
 for all $t\in[0,\kappa_0]$.
Suppose otherwise, and let $\tau'\in(0,\kappa_0]$ be the first time such that
$\bar{\omega}_x(\tau')=-\sqrt{2\nu}.$ Then $\bar{\omega}_x(t)<-\sqrt{2\nu}$ on $[0,\tau')$. From (\ref{Ktx2}), we have
$\frac{d\bar{\omega}_x}{dt}<0$, $a.\ e.$ on $[0,\tau')$, so
$\bar{\omega}_x(\tau')<\bar{\omega}_x(0)<-\sqrt{2\nu},$ a contradiction. Hence $\bar{\omega}_x(t)<-\sqrt{2\nu}$ for all $t\in[0,\kappa_0]$, and consequently $\bar{\omega}_x(t)$ is strictly decreasing.
 On the one hand, it follows from (\ref{Ktx2}) and (\ref{Ktx3}) that
\begin{eqnarray}\label{Ktx4}
\frac{d\bar{\omega}_{x}(t)}{dt}
&\leq& -\frac{\bar{\omega}^2_{x}(t)}{2}+\nu\leq (-\frac{1}{2}+\frac{\nu}{\omega^2_{0,x}})\bar{\omega}^2_{x}(t).
 \end{eqnarray}
Solving this inequality (\ref{Ktx4}) gives
\begin{eqnarray*}
\bar{\omega}_{x}(t)\leq \big(\frac{1}{\omega_{0,x}}+(\frac{1}{2}-\frac{\nu}{\omega^2_{0,x}})t\big)^{-1}.
 \end{eqnarray*}
Hence there exists $t^*:=\frac{-\frac{1}{\omega_{0,x}}}{\frac{1}{2}-\frac{\nu}{\omega^2_{0,x}}}>0$ such that $\bar{\omega}_{x}(t)$ tends to $-\infty$ as $t\rightarrow t^*$.
Therefore, $u_{x}(t,q(t,x_0))$ tends to $-\infty$ as $t\rightarrow t^*$ as a result of the boundedness of $\Omega\partial_xp\ast \rho^2.$
On the other hand, we derive from (\ref{Ktx2}) and (\ref{Ktx3}) that
\begin{eqnarray}\label{Ktx5}
\big(\frac{1}{\bar{\omega}_{x}(t)-\sqrt{2\nu}}-\frac{1}{\bar{\omega}_{x}(t)+\sqrt{2\nu}}\big)\cdot\frac{d\bar{\omega}_{x}(t)}{dt}\leq -\sqrt{2\nu}.
 \end{eqnarray}
Solving this inequality (\ref{Ktx5}) yields
\begin{eqnarray*}
t\leq \frac{1}{\sqrt{2\nu}}\big(\ln{\frac{\omega_{0,x}-\sqrt{2\nu}}{\omega_{0,x}+\sqrt{2\nu}}}-\ln{\frac{\bar{\omega}_{x}(t)-\sqrt{2\nu}}{\bar{\omega}_{x}(t)+\sqrt{2\nu}}}\big).
 \end{eqnarray*}
Taking the limit as $t\rightarrow t^*$ in the above inequality, and in view of the assumption (\ref{x0condition}), we get
\begin{eqnarray*}
T\leq \frac{1}{\sqrt{2\nu}}\ln{\frac{\omega_{0,x}-\sqrt{2\nu}}{\omega_{0,x}+\sqrt{2\nu}}}\leq \kappa_0.
 \end{eqnarray*}
This completes the proof of Theorem \ref{Th3.4}.
\end{proof}

\begin{theorem3}\label{Th3.5}
 Let
$(u_0,\rho_0-1)\in H^s\times H^{s-1}$ with $s>\frac{3}{2}$ and $T>0$ be the maximal existence time of the corresponding solution $(u,\rho)$
 to (\ref{systems2}) (guaranteed by Theorem \ref{Th3.1}). Assume that there exists $M>0$ such that $\|\rho(t,\cdot)\|_{L^\infty},\ \|\rho(t,\cdot)\|_{L^2}\le e^{2Mt}\|\rho_0\|_{H^{s-1}}$ for all $t \in [0,T)$.
For some $\kappa_0\in(0,\kappa]$ ($\kappa$ is defined in Lemma \ref{Lem3.3}) and if there exists a point $x_0\in \R$ such that
 \begin{eqnarray}\label{x0condition1}
 u_{0,x}(x_0)\leq -|u_0(x_0)|-3A_0-\sqrt{2C_1}-\sqrt{\frac{2\sqrt{2C_1}}{e^{\sqrt{2C_1}}\kappa_0-1}},
 \end{eqnarray}
 where $C_1:=\big(\frac{1}{2}+2\Omega U(\kappa_0)\big)R_{\kappa_0}^2+\frac{1}{2}U^2(\kappa_0)$
 with $U(\kappa_0)=\frac{\|\omega_0\|_{L^\infty} + A_0}
{1 - C_0\bigl(\|\omega_0\|_{L^\infty} + A_0\bigr)\kappa_0}
- \frac{1}{6}\Omega R_{\kappa_0}^2$, $R_{\kappa_0}= e^{2M{\kappa_0}}\|\rho_0\|_{H^{s-1}}$
and $A_0, C_0$ are defined as in Lemma \ref{Lem3.3} with \(\kappa\) replaced by \(\kappa_0\),
then corresponding solution $(u,\rho)$
 to (\ref{systems2}) blows up in finite time with $T\leq \kappa_0$
\end{theorem3}

\begin{proof}
Applying Theorem \ref{Th3.1} and a simple density argument, we need only to show that Theorem \ref{Th3.5} holds with some $s\geq3$.
We first introduce the following two functions
\begin{eqnarray*}
M(t):=(\omega+\omega_x)(t,q(t,x_0)) \quad \mbox{and}\quad N(t):=(\omega-\omega_x)(t,q(t,x_0)), \quad \forall t\in [0,T).
 \end{eqnarray*}
We then define the two convolution operators $p_{+}$ and $p_{-}$ as
\begin{eqnarray*}
p_{+}\ast f(x):=\frac{e^{-x}}{2}\int^x_{-\infty}e^yf(y)dy \quad \mbox{and}\quad p_{-}\ast f(x):=\frac{e^{x}}{2}\int_x^{\infty}e^{-y}f(y)dy.
 \end{eqnarray*}
Thus we obtain
\begin{eqnarray}\label{convolution}
p=p_{+}+p_{-}\quad \mbox{and}\quad \partial_xp=p_{-}-p_{+}.
 \end{eqnarray}
In view of (\ref{uut2})-(\ref{uut3}), along the trajectory of $q(t,x_0)$, we have
\begin{eqnarray}\label{At}
\frac{dM(t)}{dt}=\frac{d\omega}{dt}+\frac{d\omega_x}{dt}
&=&-3p_{-}\ast(u^2)+\frac{3}{2}u^2-u^2_x-(1-2\Omega u)p_{-}\ast(\rho^2)+(\frac{1}{2}-\Omega u)\rho^2\nonumber \\
&\leq&(u^2-u^2_x)+\frac{1}{2}u^2+2\Omega u p_{-}\ast(\rho^2)+(\frac{1}{2}-\Omega u)\rho^2,
 \end{eqnarray}
and
\begin{eqnarray}\label{Bt}
\frac{dN(t)}{dt}=\frac{d\omega}{dt}-\frac{d\omega_x}{dt}&=&3p_{+}\ast(u^2)-\frac{3}{2}u^2+u^2_x+(1-2\Omega u)p_{+}\ast(\rho^2)-(\frac{1}{2}-\Omega u)\rho^2\nonumber \\
&\geq& (-u^2+u^2_x)-\frac{1}{2}u^2-2\Omega u p_{+}\ast(\rho^2)-(\frac{1}{2}-\Omega u)\rho^2.
 \end{eqnarray}
Applying the a priori estimate (\ref{estimateu1}), Young's inequality and the assumptions on $\rho$, we get the following estimates for all $t\in [0,\kappa_0]\cap [0,T)$
 \begin{eqnarray}\label{ABt}
&&|\pm\frac{1}{2}u^2\pm2\Omega u p_{\mp}\ast(\rho^2)\pm(\frac{1}{2}-\Omega u)\rho^2|\nonumber \\
&\leq & \big(\frac{1}{2}+2\Omega U(\kappa_0)\big)R_{\kappa_0}^2+\frac{1}{2}U^2(\kappa_0):=C_1.
 \end{eqnarray}
 where $U(\kappa_0)=\frac{\|\omega_0\|_{L^\infty} + A_0}
{1 - C_0\bigl(\|\omega_0\|_{L^\infty} + A_0\bigr)\kappa_0}
- \frac{1}{6}\Omega R_{\kappa_0}^2$, $R_{\kappa_0}= e^{2M{\kappa_0}}\|\rho_0\|_{H^{s-1}}$
and $A_0, C_0$ are defined as in Lemma \ref{Lem3.3} with \(\kappa\) replaced by \(\kappa_0\).
Substituting (\ref{ABt}) into (\ref{At})-(\ref{Bt}), we deduce that
 \begin{eqnarray}\label{AtB}
\frac{dM(t)}{dt}\leq u^2-u^2_x+C_1 \quad \mbox{and}\quad \frac{dN(t)}{dt}\geq -u^2+u^2_x-C_1.
 \end{eqnarray}
In addition, it follows from the definition of $\omega$, Young's inequality and (\ref{convolution}) that
\begin{eqnarray}\label{AtB1}
u+u_x \leq M(t) \leq u+u_x +\Omega R_{\kappa_0}^2 =u+u_x +\frac{3}{2}A_0,
 \end{eqnarray}
and
\begin{eqnarray}\label{AtB2}
 u-u_x \leq N(t) \leq u-u_x +\Omega R_{\kappa_0}^2=u-u_x +\frac{3}{2}A_0.
 \end{eqnarray}

By our assumption (\ref{x0condition1}) on the initial data, we have
\begin{eqnarray*}
u^2_{0,x}(x_0)-u_0^2(x_0)-C_1>0,
 \end{eqnarray*}
which implies, from (\ref{AtB}), that
 \begin{eqnarray}\label{AtB0}
\frac{d}{dt}M(0)<0 \quad \mbox{and}\quad \frac{d}{dt}N(0)>0.
 \end{eqnarray}
Moreover, from the assumption (\ref{x0condition1}) and (\ref{AtB1})-(\ref{AtB2}), we find that
\begin{eqnarray}\label{AtB00}
 M(0)<-\sqrt{2C_1}-\frac{3}{2}A_0<0 \quad \mbox{and}\quad N(0)>\sqrt{2C_1}+3A_0>0.
 \end{eqnarray}
Thus over the time of existence we claim that
\begin{eqnarray}\label{AtB01}
\frac{dM(t)}{dt}<0  \quad \mbox{and}\quad  \frac{dN(t)}{dt}>0.
 \end{eqnarray}
Assume the above monotonicity (\ref{AtB01}) fails at some time. From (\ref{AtB0}), by continuity, let $[0,\kappa_0]\cap [0,T)\ni t_0>0$
be the first time such that
\begin{eqnarray*}
\frac{d}{dt}M(t_0)=0 \quad \mbox{or}\quad \frac{d}{dt}N(t_0)=0.
 \end{eqnarray*}
Since the two cases are symmetric, we only prove the case $\frac{d}{dt}M(t_0)=0$.
On $[0,t_0)$, $M(t)$ is strictly decreasing and $N(t)$ strictly increasing, hence from (\ref{AtB00}), we get
\begin{eqnarray}\label{AtB02}
M(t_0)<M(0)<-\sqrt{2C_1}-\frac{3}{2}A_0<0 \quad \mbox{and}\quad N(t_0)>N(0)>\sqrt{2C_1}+3A_0>0.
 \end{eqnarray}
Combining (\ref{AtB1})-(\ref{AtB2}) with (\ref{AtB02}), we obtain
\begin{eqnarray*}
(u+u_x)(t_0) \leq M(t_0) <-\sqrt{2C_1}-\frac{3}{2}A_0 \quad \mbox{and}\quad  (u-u_x)(t_0)\geq N(t_0)-\frac{3}{2}A_0>\sqrt{2C_1}+\frac{3}{2}A_0,
 \end{eqnarray*}
which, together with the first inequality in (\ref{AtB}), gives
\begin{eqnarray*}
\frac{d}{dt}M(t_0)\leq (u^2-u^2_x)(t_0)+C_1<-(\sqrt{2C_1}+\frac{3}{2}A_0)^2+C_1<0.
 \end{eqnarray*}
This contradicts the assumption $\frac{d}{dt}M(t_0)=0$. Hence, the claim (\ref{AtB01}) holds on $[0,\kappa_0]\cap [0,T)$.
Therefore, from the assumption (\ref{x0condition1}) and (\ref{AtB1})-(\ref{AtB2}), we know that
\begin{eqnarray}\label{AtB03}
M(t)<M(0)\leq u_0+u_{0,x}+\frac{3}{2}A_0
<-\frac{3}{2}A_0-\sqrt{2C_1}-\sqrt{\frac{2\sqrt{2C_1}}{e^{\sqrt{2C_1}}\kappa_0-1}}<0,
\end{eqnarray}
and
\begin{eqnarray}\label{AtB033}
N(t)>N(0)\geq u_0-u_{0,x}> 3A_0+\sqrt{2C_1}+\sqrt{\frac{2\sqrt{2C_1}}{e^{\sqrt{2C_1}}\kappa_0-1}}>0,
 \end{eqnarray}
which allows us to define the function
\begin{eqnarray*}
h(t)=\sqrt{-M(t)N(t)}>0.
 \end{eqnarray*}
Computing the derivative of $h$, next applying the differential inequalities (\ref{AtB}) and (\ref{AtB03})-(\ref{AtB033}), and the inequality $\frac{N(t)-M(t)}{2}\geq h(t)$,
we have
\begin{eqnarray}\label{AtB04}
\frac{d h(t)}{dt}=-\frac{\frac{dM(t)}{dt}
N+M\frac{dN(t)}{dt}}{2\sqrt{-M(t)N(t)}}\geq (u_x^2-u^2-C_1)\frac{N(t)-M(t)}{2\sqrt{-M(t)N(t)}}\geq u_x^2-u^2-C_1.
 \end{eqnarray}
In view of (\ref{AtB1})-(\ref{AtB2}) and (\ref{AtB03}), we obtain
\begin{eqnarray*}
0<- M(t) \leq -(u+u_x)\quad \mbox{and}\quad 0<N(t)\leq u-u_x +\frac{3}{2}A_0 \leq 2(u-u_x),
 \end{eqnarray*}
which yields
\begin{eqnarray}\label{AtB05}
h^2(t)=-M(t)N(t)\leq 2(u^2_x-u^2).
 \end{eqnarray}
It then follows from (\ref{AtB04})-(\ref{AtB05}) that
\begin{eqnarray}\label{AtB06}
\frac{d h(t)}{dt}\geq \frac{1}{2}h^2-C_1.
 \end{eqnarray}
But at initial time, from (\ref{AtB03})-(\ref{AtB033}), we see that
\begin{eqnarray}\label{AtB07}
h^2(0)=-M(0)N(0)>\big(\frac{3A_0}{2}+\sqrt{2C_1}+\sqrt{\frac{2\sqrt{2C_1}}{e^{\sqrt{2C_1}}\kappa_0-1}}
\big)^2
\end{eqnarray}
Applying a similar argument used for (\ref{Ktx3}) in Theorem \ref{Th3.4}, one infers that
there exists $t^*:=\frac{\frac{1}{h(0)}}{\frac{1}{2}-\frac{C_1}{h^2(0)}}>0$
such that $h(t)$ tends to $+\infty$ as $t\rightarrow t^*$. Note that
\begin{eqnarray*}
h(t)=\sqrt{\omega^2_x-\omega^2}\leq -\omega_x=-u_x-\Omega\,\partial_xp\ast\rho^2\leq-u_x+\frac{1}{2}\Omega R^2_{\kappa_0}.
 \end{eqnarray*}
Therefore, $u_{x}(t,q(t,x_0))$ tends to $-\infty$ as $t\rightarrow t^*$.
On the other hand, it follows from (\ref{AtB06}) and (\ref{AtB07}) that
\begin{eqnarray*}
\big(\frac{1}{h(t)-\sqrt{2C_1}}-\frac{1}{h(t)+\sqrt{2C_1}}\big)\cdot\frac{dh(t)}{dt}\geq \sqrt{2C_1}.
 \end{eqnarray*}
Solving the above inequality gives
\begin{eqnarray*}
t\leq \frac{1}{\sqrt{2C_1}}\big(\ln{\frac{h(0)+\sqrt{2C_1}}{h(0)-\sqrt{2C_1}}}-\ln{\frac{h(t)+\sqrt{2C_1}}{h(t)-\sqrt{2C_1}}}\big).
 \end{eqnarray*}
Taking the limit as $t\rightarrow t^*$ in the above inequality, and by the condition (\ref{AtB07}), we get
\begin{eqnarray*}
T\leq \frac{1}{\sqrt{2C_1}}\big(\ln{\frac{h(0)+\sqrt{2C_1}}{h(0)-\sqrt{2C_1}}}\big)\leq \kappa_0.
 \end{eqnarray*}
This completes the proof of Theorem \ref{Th3.5}.
 \end{proof}

\begin{remark3}\label{Re3.1}
In Theorems \ref{Th3.4}-\ref{Th3.5}, if the initial data satisfy that at the same point $x_0\in \R$, we have
\begin{eqnarray*}
\rho_0(x_0)=0
 \end{eqnarray*}
and $u_{0,x}(x_0)$ satisfies assumption (\ref{x0condition}) or (\ref{x0condition1}), then the a priori bound
on $\|\rho(t,\cdot)\|_{L^\infty}$ (i.e. $\|\rho(t,\cdot)\|_{L^\infty}\le e^{2Mt}\|\rho_0\|_{H^{s-1}}$ for all $t\in[0,T)$) can be completely removed.
Indeed, from the second equation of system (\ref{systems2}) along the characteristic $q(t,x)$ defined in (\ref{flow1}), we obtain
\begin{eqnarray*}
\frac{d\rho(t,q(t,x_0))}{dt}=-2(\rho u_x)(t,q(t,x_0)),
 \end{eqnarray*}
and since $\rho(0,q(0,x_0))=\rho_0(x_0)=0$, it follows that
\begin{eqnarray*}
\rho(t,q(t,x_0))\equiv 0, \qquad \forall t\in[0,T).
 \end{eqnarray*}
Consequently, the term $\left(\frac{1}{2}-\Omega u\right)\rho^2$ vanishes identically, and its estimate no longer requires any bound on $\|\rho\|_{L^\infty}$. All other terms, such as $(\partial_x p*\rho^2)^2$ and $ u\,p*\rho^2$, depend only on $\|\rho\|_{L^2}$ and $\|u\|_{L^\infty}$, and the latter is controlled by Lemma \ref{Lem3.3} under the sole assumption of a bound on $\|\rho\|_{L^2}$.
Thus, under the condition that both $\rho_0(x_0)=0$ and (\ref{x0condition}) or (\ref{x0condition1}) hold at the same point $x_0$, Theorems \ref{Th3.4}-\ref{Th3.5} remain valid without the extra $\|\rho(t,\cdot)\|_{L^\infty}$ hypothesis.
\end{remark3}

\begin{remark3}\label{Re3.2}
When the parameters in the R-$b$-family system (\ref{systemsoriginal}) are set to $b=3,\ \mu=0,\ A=0,\ \sigma=1$, it reduces to
\begin{equation}\label{systems6}
\left\{
 \begin{aligned}
&m_t+um_x+3u_xm+\rho\rho_x-2\Omega \rho(\rho u)_x=0,\quad m=u-u_{xx},\\
&\rho_t+(\rho u)_x=0.
\end{aligned}
\right.
\end{equation}
This is a special case of the system (\ref{systemsoriginal}) with the similar nonlinear structure as the R2DP system \eqref{systems2}.
Hence, by introducing the same auxiliary variable $\omega= u + \Omega\,p\ast\rho^2,\ p(x)=\frac12 e^{-|x|},$ we obtain
 \begin{eqnarray*}
\omega_t+u\omega_x
=\Omega \omega \partial_xp\ast (\rho^2)-\Omega^2(p\ast\rho^2)(\partial_xp\ast\rho^2)-\partial_xp\ast (\frac{3}{2}u^2+\frac{1}{2}\rho^2),
 \end{eqnarray*}
 and
  \begin{eqnarray*}
\omega_{tx}+u\omega_{xx}
&=&-\big(\omega_x-\Omega\,\partial_xp\ast\rho^2\big)^2+(\frac{1}{2}-\Omega u)\rho^2+\Omega u p\ast\rho^2
+\frac{3}{2}u^2-p\ast(\frac{3}{2}u^2+\frac{1}{2}\rho^2).
 \end{eqnarray*}
 Hence although the two systems \eqref{systems2} and \eqref{systems6} differ in their explicit form, the coupling term in both cases can be rewritten as
 $-\rho\rho_t$ using the respective $\rho$-equation, so the resulting $\omega$-equations are identical. Consequently,
through the exact same procedure as in Sections \ref{Sec3.1}-\ref{Sec3.2}, the blow-up scenario (Theorems \ref{Th3.2}-\ref{Th3.3}
), and the blow-up criterion (Theorems \ref{Th3.4}-\ref{Th3.5}) all hold for the system \eqref{systems6} accordingly.
\end{remark3}

\begin{remark3}\label{Re3.3}
All results established in Sections \ref{Sec3.1}-\ref{Sec3.2}, and Remarks \ref{Re3.1}-\ref{Re3.2} for the R2DP system \eqref{systems2} and the system \eqref{systems6} remain valid in the periodic setting, where the spatial variable belongs to the unit circle \(\mathbb S=\mathbb R/\mathbb Z\).
In this case, the kernel of the operator \((1-\partial_x^2)^{-1}\) is given by
\[
G(x)=\frac{\cosh\left(x-[x]-\frac12\right)}{2\sinh\left(\frac12\right)},\qquad x\in\mathbb S.
\]
\end{remark3}

\begin{remark3}\label{Re3.4}
When $\Omega=0$, the R2DP system (\ref{systems2}) reduces to 2DP system.
In that case, the 2DP system admits an a priori estimate for $\|u\|_{L^\infty}$ on the whole
maximal existence interval $[0,T)$ (see Lemma 4.3 in \cite{YY}).
In contrast, the $L^\infty$-estimate for $u$ established in our Lemma \ref{Lem3.3}
is only available on the local time interval $[0,\kappa]$.
Consequently, when $\Omega=0$ the validity of the $L^\infty$-bound on the whole lifespan simplifies the application of Theorem \ref{Th3.4},
and our result recovers the corresponding wave-breaking criterion for the 2DP system (see Theorem 5.1 in \cite{YY}).
Furthermore, if we additionally set $\rho\equiv0$, the system further reduces to the DP equation, and Theorem \ref{Th3.4}
coincides with the blow-up criterion of Theorem 3.1 in \cite{ZY}.
Moreover, to the best of our knowledge, Theorem \ref{Th3.5} provides new wave-breaking criteria even for the 2DP system and DP equation.
\end{remark3}

\noindent\textbf{Acknowledgments} \
The authors thank the anonymous referee for helpful suggestions and comments.
The work is supported by National Natural Science Foundation of China under Grant 12001528.\\

\noindent\textbf{Conflicts of Interest} \
The authors declare no conflicts of interest.\\

\noindent\textbf{Data Availability Statement} \
Data sharing not applicable to this article as no datasets were generated or analyzed during the current study.\\

\end{document}